\documentclass[letterpaper,12pt]{amsart}

\usepackage{amsmath,amsthm,amssymb,epic,epsfig,eepic,mathtools,setspace,dsfont,hyperref,color,graphicx}
\renewcommand{\thispagestyle}[1]{} 
\renewcommand{\qed}{$\blacksquare$}
\newcommand{\dia}{$\diamond$}

\newtheorem{thm}{Theorem}[section]

\newtheorem{cor}[thm]{Corollary}
\newtheorem{lemma}[thm]{Lemma}
\newtheorem{prop}[thm]{Proposition}

\newtheorem{dfn}[thm]{Definition} 
\newtheorem{rem}[thm]{Remark} 

\newtheorem{ex}[thm]{Example}

\newcommand{\thth}{^{\text{\underline{th}}}}

\newcommand{\conv}{\mathrm{Conv}}
\newcommand{\rank}{\mathrm{Rank}}

\newcommand{\C}{\mathbb{C}}

\newcommand{\Cn}{\C^n}

\newcommand{\Pro}{\mathbb{P}}
\newcommand{\R}{\mathbb{R}}
\newcommand{\Rs}{\R^*}
\newcommand{\Rn}{\R^n}
\newcommand{\Z}{\mathbb{Z}}

\newcommand{\cA}{\mathcal{A}}
\newcommand{\hA}{\widehat{\cA}}

\newcommand{\sign}{\mathrm{sign}}
\newcommand{\Log}{\mathrm{Log}}

\newcommand{\cC}{\mathcal{C}}

\newcommand{\cL}{\mathcal{L}}

\makeatletter
\renewcommand*\env@matrix[1][c]{\hskip -\arraycolsep
  \let\@ifnextchar\new@ifnextchar
  \array{*\c@MaxMatrixCols #1}}
\makeatother

\begin{document}

\title[Subexponential Fewnomial Hypersurface Bounds]{\mbox{}
\vspace{-1in}\\ 
New Subexponential Fewnomial Hypersurface Bounds} 

\author{Jens Forsg\aa{}rd} 
\address{Department of Mathematics,
Texas A\&M University TAMU 3368,
College Station, Texas \ 77843-3368, USA.}
\email{jensf@math.tamu.edu }
\author{Mounir Nisse} 
\address{1415, School of Mathematics, Korea Institute for Advanced Study 
(KIAS), 85 Hoegiro, Dongdaemun-gu, Seoul, 02455, Republic of Korea.}
\email{mounir.nisse@gmail.com}
\author{J.\ Maurice Rojas}
\address{Department of Mathematics,
Texas A\&M University TAMU 3368,
College Station, Texas \ 77843-3368, USA.}
\email{rojas@math.tamu.edu}
\thanks{Partially supported by NSF grant CCF-1409020 and MSRI. } 

\begin{abstract} Suppose $c_1,\ldots,c_{n+k}$ are real 
numbers, $\{a_1,\ldots,a_{n+k}\}\!\subset\!\Rn$ is a set of points 
not all lying in the same affine hyperplane, $y\!\in\!\Rn$, 
$a_j\cdot y$ denotes the standard real inner product of $a_j$ and 
$y$, and we set $g(y)\!:=\!\sum^{n+k}_{j=1} c_j e^{a_j\cdot y}$. We prove 
that, for generic $c_j$, the number of connected components of the real 
zero set of $g$ is $O\!\left(n^2+\sqrt{2}^{k^2}(n+2)^{k-2}\right)$. 
The best previous upper bounds, when restricted to the special case 
$k\!=\!3$ and counting just non-compact components, were already exponential 
in $n$.  
\end{abstract}  

\maketitle 

\vspace{-.3in} 
\section{Introduction} 
Estimating the number of connected components of the real zero set 
of a system of polynomial equations is a fundamental problem occuring 
in numerous applications. For instance, in robotics \cite{wms,chase}, chemical 
reaction networks \cite{shiu}, economic modelling \cite{mclennan}, and 
complexity theory \cite{koiran}, information on the topology of the underlying 
zero set is sometimes 
at least as important as numerically approximating solutions.  
We derive topological bounds in the broader context of real exponential 
sums, significantly sharpening older bounds from fewnomial theory 
\cite{khovanski,bihansottilestrat}.  
\begin{dfn} 
\label{dfn:a} 
For any field $K$ we let $K^*\!:=\!K\setminus\{0\}$. Let 
$\cA\!\in\!\R^{n\times (n+k)}$ have $j\thth$ column $a_j$   
and let $c_1,\ldots,c_{n+k}\!\in\!\Rs$. We then call 
$g(y)\!:=\!\sum^t_{j=1} c_j e^{a_j\cdot y}$ a {\em (real) 
$n$-variate exponential $(n+k)$-sum}, and call 
$\cA$ the {\em spectrum} of $g$. We also let $c_g\!:=\!(c_1,\ldots,c_{n+k})$. 
Finally, for any function $h: \Cn\longrightarrow \R$, we let 
$Z_\C(h)$, $Z_\R(h)$, and $Z_+(h)$ respectively denote the zeroes of 
$h$ in $\Cn$, $\Rn$, and $\Rn_+$ (the positive orthant). \dia  
\end{dfn}
 
\noindent 
Note that when $\cA\!\in\!\Z^{n\times (n+k)}$ 
there is an obvious $f\!\in\!\R\!\left[x^{\pm 1}_1,\ldots,x^{\pm 1}_n\right]$, 
with exactly $n+k$ monomial terms, such that $g(y)\!=\!f\!\left(e^{y_1},
\ldots,e^{y_n}\right)$ identically, and the zero sets $Z_\R(g)$ and $Z_+(f)$ 
have the same number of connected components. In this sense, among many 
others, real exponential sums generalize real polynomials. 

We say a condition involving a tuple of real parameters                 
$(z_1,\ldots,z_N)$ holds {\em generically} if and only if the set of
choices of $(z_1,\ldots,z_N)$ making the condition true is dense and open
in $\R^N$. For instance, it is easy to show that for generic 
$\cA\!\in\!\R^{n\times (n+k)}$ (with $k\!\geq\!1$) we have  
that $\{a_1,\ldots,a_{n+k}\}$ do not all lie in the same affine 
hyperplane. 
\begin{thm}
\label{thm:big} 
Suppose $g$ is an $n$-variate $(n+k)$-sum with spectrum $\cA$ and  
$\{a_1,\ldots,a_{n+k}\}$ do not all lie in the same affine hyperplane. 
Then, for generic $c_g$, $Z_\R(g)$ has no more than 
$\displaystyle{\frac{(n+k)(n+k-1)}{2}+\left\lfloor 
\sqrt{2}^{(k-2)(k-3)}(n+2)^{k-2}\right\rfloor}$ 
connected components. Furthermore, for $k\!=\!3$, a sharper upper bound of 
$\frac{(n+3)(n+2)}{2}+\left\lfloor \frac{n+5}{2}\right\rfloor$ holds. 
\end{thm} 

\noindent 
We prove Theorem \ref{thm:big} in Section \ref{sub:proof} below. 
The best previous upper bound on the number of connected components,  
\cite[Thm.\ 1]{bihansottilestrat}, came 
from a larger topological invariant: the sum of the {\em Betti numbers} of the 
underlying zero set. (See also \cite{basu} for an important precursor in the 
semi-algebraic setting.)  
Our bound is polynomial in $n$ for any fixed $k$, while the bound 
from \cite[Thm.\ 1]{bihansottilestrat} is exponential in each of $n$ 
and $k$. For $k\!\in\!\{1,2\}$ respective optimal 
upper bounds of $1$ and $2$ are already known (see, e.g., 
\cite{brs,bihan,reutopo}). 

\section{Background} 
\label{sec:back}  
A central tool behind the proof of Theorem \ref{thm:big} 
is an extension of Gelfand, Kapranov, and Zelevinsky's theory of 
$\cA$-discriminants \cite{gkz94} to exponential sums. This generalization was 
first developed in \cite{rojasrusek}.  
\begin{dfn} 
For any $\cA\!\in\!\R^{n\times (n+k)}$ we define the {\em generalized 
$\cA$-discriminant variety}, $\Xi_\cA$, 
to be the Euclidean closure of the set
of all $[c_1:\cdots:c_{n+k}]\!\in\!\Pro^{n+k-1}_\C$ such that
$\sum^{n+k}_{j=1} c_je^{a_j\cdot z}$ has a degenerate root in $\Cn$. 
Also, we call $\cA$ {\em non-defective} if and only if $\Xi_\cA$ has 
codimension $1$ in $\Pro^{n+k-1}_\C$. \dia 
\end{dfn} 
\begin{dfn}
\label{dfn:iso}
Given any two subsets $X,Y\!\subseteq\!\Rn$, an {\em isotopy from $X$ to $Y$
(ambient in $\Rn$)} is a continuous map $I : [0,1]\times \Rn \longrightarrow 
\Rn$ satisfying (1) $I(t,\cdot)$ is a homeomorphism for all
$t\!\in\![0,1]$, (2) $I(0,x)\!=\!x$ for all $x\!\in\!\Rn$, and
(3) $I(1,X)\!=\!Y$. \dia
\end{dfn}

\noindent
It is easily checked that an isotopy from $X$ to $Y$ implies an isotopy
from $Y$ to $X$ as well.
So isotopy is in fact an equivalence relation and it makes sense to speak of
{\em isotopy type}.

The real part of $\Xi_\cA$ (along with some additional pieces: see 
Theorems \ref{thm:morse} and \ref{thm:morse2} below) partitions the coefficient 
space of $g$ into regions where $Z_\R(g)$ is smooth and the isotopy type of 
$Z_\R(g)$ is constant. Moreover, since scaling variables and coefficient 
vectors does not affect the presence of singularities in $Z_\R(g)$, the variety 
$\Xi_\cA$ has certain homogeneities. 
As we'll see below, these homogeneities can be quotiented out to better study 
regions of the coefficient space where $Z_\R(g)$ is smooth and has constant 
isotopy type. For any $S\!\subseteq\!\C^N$ we let $\overline{S}$ 
denote the Euclidean closure of $S$. 
\begin{dfn}
\label{dfn:first} 
For any $\cA\!\in\!\R^{n\times (n+k)}$ let
$\hA\!\in\!\R^{(n+1)\times (n+k)}$ denote the matrix with first row\linebreak 
\scalebox{.95}[1]{$[1,\ldots,1]$ and bottom $n$ rows forming $\cA$, and set  
$d(\cA)\!:=\!\rank \hA-1$. Let $B\!\in\!\R^{(n+k)\times (n+k-d(\cA)-1)}$}
\linebreak 
be any matrix whose columns form a basis for the right nullspace of $\hA$. Let 
$\beta_i$ denote the $i\thth$ row of $B$, let $(\cdot)^\top$ denote matrix 
transpose, and for any $z\!=\!(z_1,\ldots,z_N)$ let $\Log|z|\!:=\!
(\log|z_1|,\ldots,\log|z_N|)$. When $\cA$ is non-defective we then set 
$\lambda\!:=\!(\lambda_1,\ldots,\lambda_{n+k-d(\cA)-1})$ 
and define the (projective) hyperplane arrangement\\  
\mbox{}\hfill $H_\cA\!:=\!\left\{[\lambda] \; \left| 
\; \; \lambda\cdot \beta_i\!=\!0\text{ for some (nonzero) row } \beta_i 
\text{ of } B\right. \right\}\subset\Pro^{n+k-d(\cA)-2}_\C$.\hfill\mbox{}\\ 
Finally, we define $\xi_{\cA,B} : 
\left(\left. \Pro^{n+k-d(\cA)-2}_\C\right\backslash \!H_\cA\right) \
\longrightarrow \R^{n+k-d(\cA)-1}$ by 
$\xi_{\cA,B}([\lambda])\!:=\!\left(\Log\left|\lambda B^\top
\right|\right)B$. (So $\xi_{\cA,B}$ is defined by multiplying a 
row vector by a matrix.) 
We then call $\Gamma(\cA,B)\!:=\!\overline{
\xi_{\cA,B}\!\left(\left. \Pro^{n+k-d(\cA)-2}_\R\right
\backslash\!H_\cA\right)}$ a {\em reduced discriminant contour}. 
\dia 
\end{dfn}

For any subset $S\!\subseteq\!\Rn$, we let $\conv S$ denote the smallest
convex set containing $S$. It is easily checked that $\dim \conv\{
a_1,\ldots,a_{n+k}\}\!=\!d(\cA)$ and thus, for generic $\cA$, we 
have $d(\cA)\!=\!n$. However, we will need to consider arbitrary 
$d(\cA)$ in order to more easily describe our approach to counting isotopy 
types. Let us call $\cA$ {\em pyramidal} if and only if $\cA$ has a 
column $a_j$ such that $\{a_1,\ldots,a_{j-1},a_{j+1},\ldots,a_{n+k}\}$
lies in a $(d(\cA)-1)$-dimensional affine subspace. The following proposition, 
on certain exceptional spectra $\cA$, will prove useful later on. 
\begin{prop}
\label{prop:pyra}
Following the preceding notation, $\cA$ is pyramidal  
if and only if $B$ has a zero row. In particular, 
$\cA$ non-defective implies that $\cA$ is not pyramidal. \qed 
\end{prop}
\begin{rem} 
When $\cA\!\in\!\Z^{n\times (n+k)}$ and $\cA$ is non-defective 
it follows easily from the development of \cite{rojasrusek} that 
$\overline{\xi_{\cA,B}\!\left(\left. 
\Pro^{n+k-d(\cA)-2}_\C\right\backslash\! H_\cA\right)}$ is in fact a linear 
section of the {\em amoeba} of the classical {\em $\cA$-discriminant 
polynomial} $\Delta_\cA$. $\xi_{\cA,B}$ is thus 
a generalization of the (logarithmic) {\em Horn-Kapranov Uniformization}  
(see \cite{kapranov,gkz94}). See also \cite{passare} for further 
background on $\cA$-discriminant contours in the special case 
$\cA\!\in\!\Z^{n\times (n+k)}$. \dia 
\end{rem} 
\begin{thm}
\label{thm:getchambers} \cite{rojasrusek} 
If $\cA$ is non-defective then $\Gamma(\cA,B)$ is a finite union of 
codimension $1$ smooth semi-analytic subsets of $\R^{n+k-d(\cA)-1}$. 
Furthermore, there is a codimension-$2$ semi-analytic set 
$Y\!\subset\!\R^{n+k-d(\cA)-1}$ such that $\Gamma(\cA,B) \cup 
Y\!=\!\left(\Log\left|\Xi_\cA \cap \Pro^{n+k-1}_\R\right|\right)B$. \qed 
\end{thm}  

\noindent 
\vbox{
\noindent 
\begin{minipage}[t]{.6 \textwidth}
\vspace{0pt}
\begin{ex} \label{ex:penta} 
When $\cA:=$\scalebox{1}[.6]{$\begin{bmatrix}0 & 1 & 0 & 4 & 1\\ 
0 & 0 & 1 & 1 & 4\end{bmatrix}$} we are in essence considering the 
family of exponential sums $g(y)\!:=\!f\!\left(e^{y_1},e^{y_2}\right)$ 
where $f(x)\!=\!c_1+c_2x_1+c_3x_2+c_4x^4_1 x_2+c_5 x_1 x^4_2$. A suitable 
$B$ (among many others) with columns defining a basis for the right nullspace 
of $\hA$ is then $B\!\approx\; $\scalebox{.6}[.6]
{$\begin{bmatrix}[r] 0.5079 & -0.8069 & 0.1721 & 0.2267 & -0.0997\\ 
0.5420 & 0.1199 & -0.7974 & -0.0851 & 0.2206 \end{bmatrix}^\top$}, and 
the corresponding reduced contour $\Gamma(\cA,B)$, intersected with $[-4,4]$,  
is drawn to the right. \dia 
\end{ex} 
In what follows, we set\\ 
\mbox{}\hfill $\sign(c_g)\!:=\!(\sign(c_1),\ldots,
\sign(c_{n+k}))\!\in\!\{\pm 1\}^{n+k}$.\hfill\mbox{}  
\end{minipage} \hspace{.3cm} 
\begin{minipage}[t]{.2 \textwidth}  
\raisebox{-5.5cm}{\epsfig{file=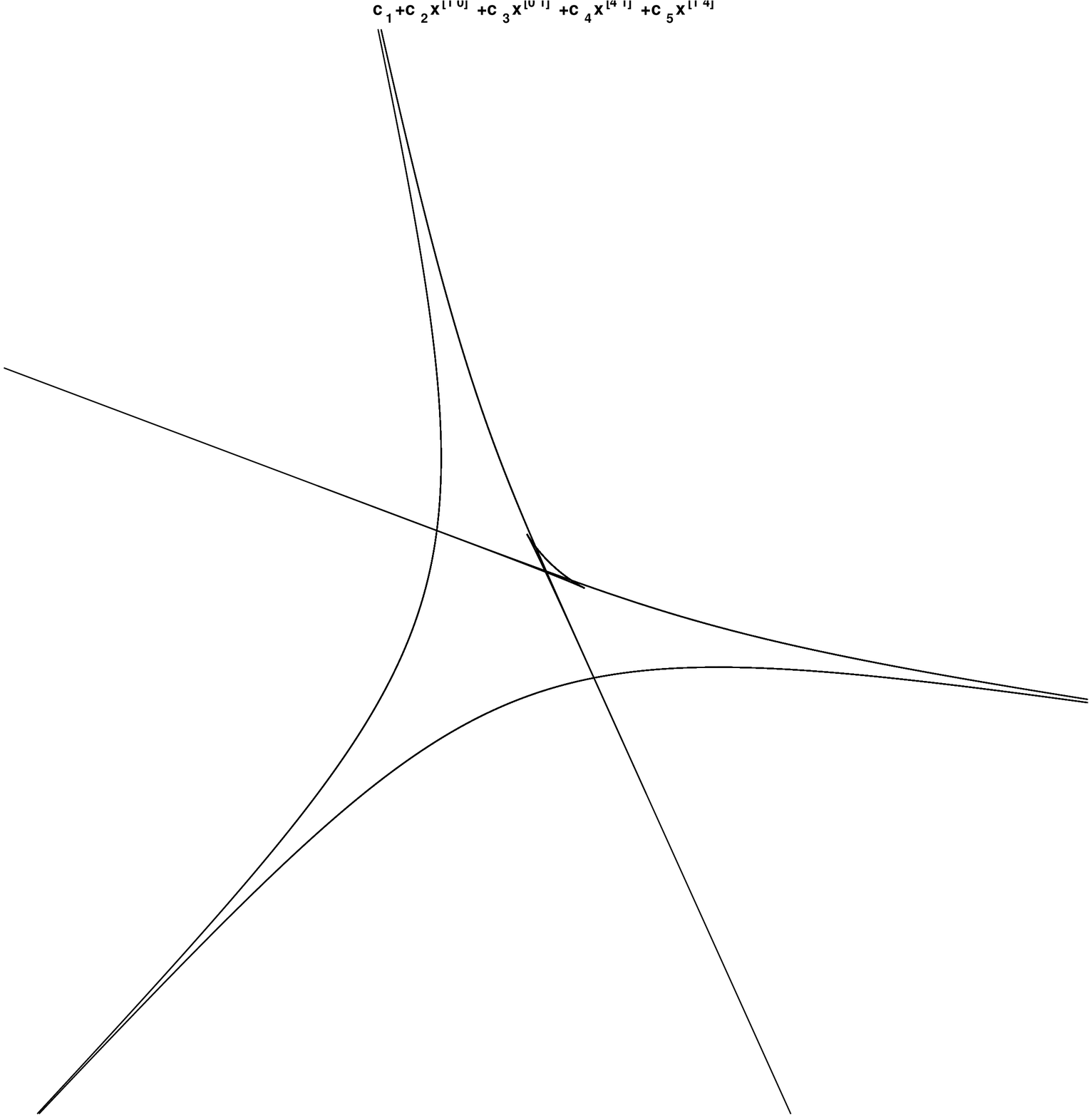,height=2.3in,clip=}}  
\end{minipage}}   
\begin{dfn}
\label{dfn:chamber}
Suppose $\cA\!\in\!\R^{n\times (n+k)}$ is non-defective and   
$\sigma\!=\!(\sigma_1,\ldots,\sigma_{n+k})\!\in\!\{
\pm 1\}^{n+k}$. We then call\\ 
\mbox{}\hfill $\Gamma_\sigma(\cA,B)\!:=\! 
\overline{\left\{\left.\xi_{\cA,B}([\lambda]) \; \right| \; 
\mathrm{sign}\left(\Log\left|\lambda B^\top\right|\right) 
\!=\!\pm \sigma  \ , \ [\lambda]\!\in\!
\Pro^{n+k-d(\cA)-2}_\R\setminus H_\cA\right\}}\!\subset\!\R^{n+k-d(\cA)-1}$
\hfill\mbox{}\\ 
a {\em signed reduced contour}, and we call any connected component $\cC$ of 
$\R^{n+k-d(\cA)-1}\setminus\Gamma_\sigma(\cA,B)$ a {\em reduced signed 
chamber}. We also call $\cC$ an {\em outer} or {\em inner chamber}, according 
as $\cC$ is unbounded or bounded. 
\dia 
\end{dfn}
\begin{ex} 
Continuing Example \ref{ex:penta}, there are $16$ possible 
choices for $\sigma$, if we identify\linebreak 
\scalebox{.96}[1]{sign sequences with their negatives. Among 
these choices, there are $11$ $\sigma$ yielding 
$\Gamma_\sigma(\cA,B)\!=\!\emptyset$.}\linebreak  
The remaining choices, along with their respective $\Gamma_\sigma(\cA,B)$ 
are drawn below. \dia 
\end{ex} 

\noindent 
\scalebox{.92}[.92]{
\begin{picture}(200,90)(12,-17)
\put(10,-15){
\fbox{\epsfig{file=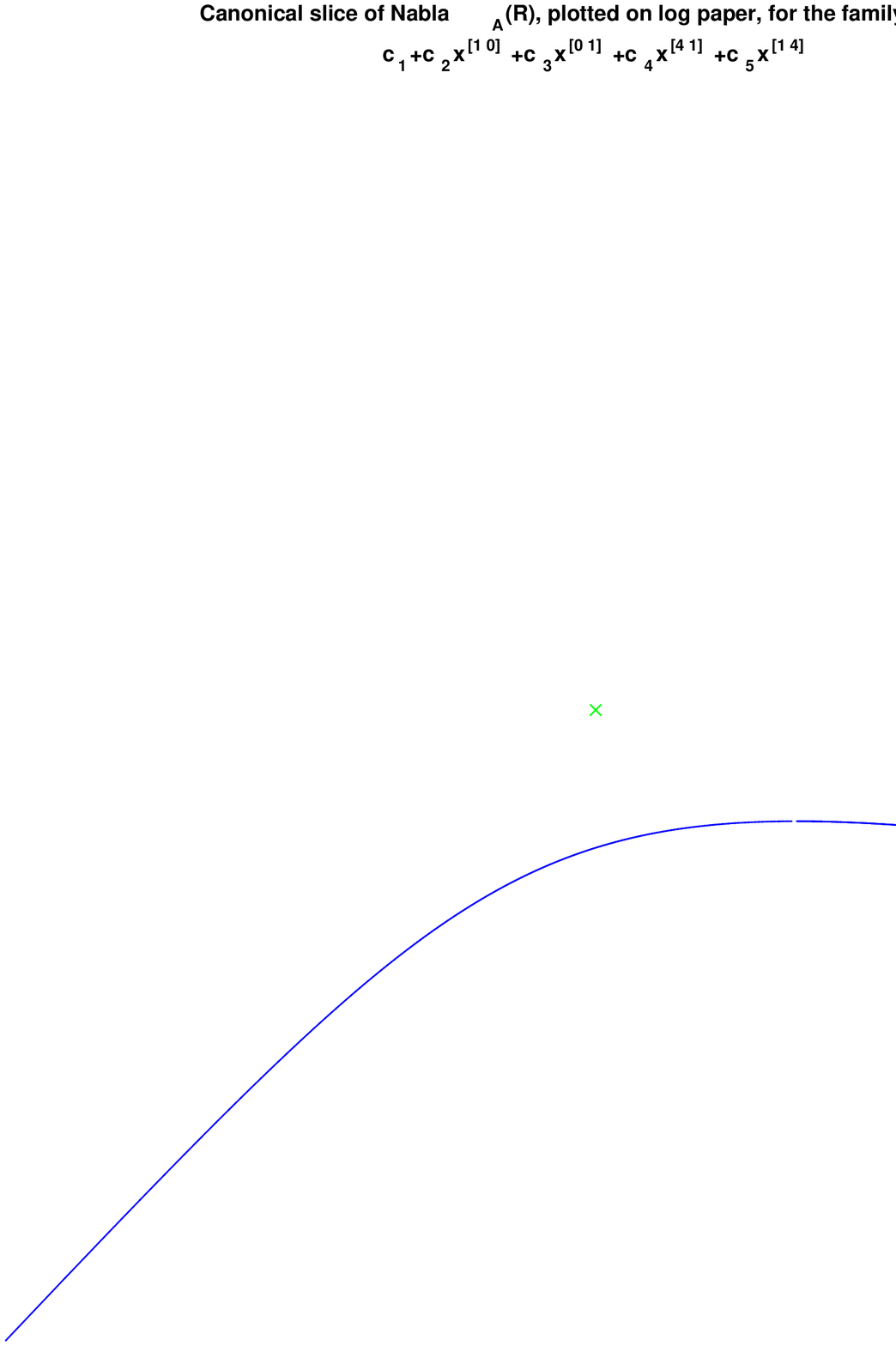,height=1.2in,clip=}}
\fbox{\epsfig{file=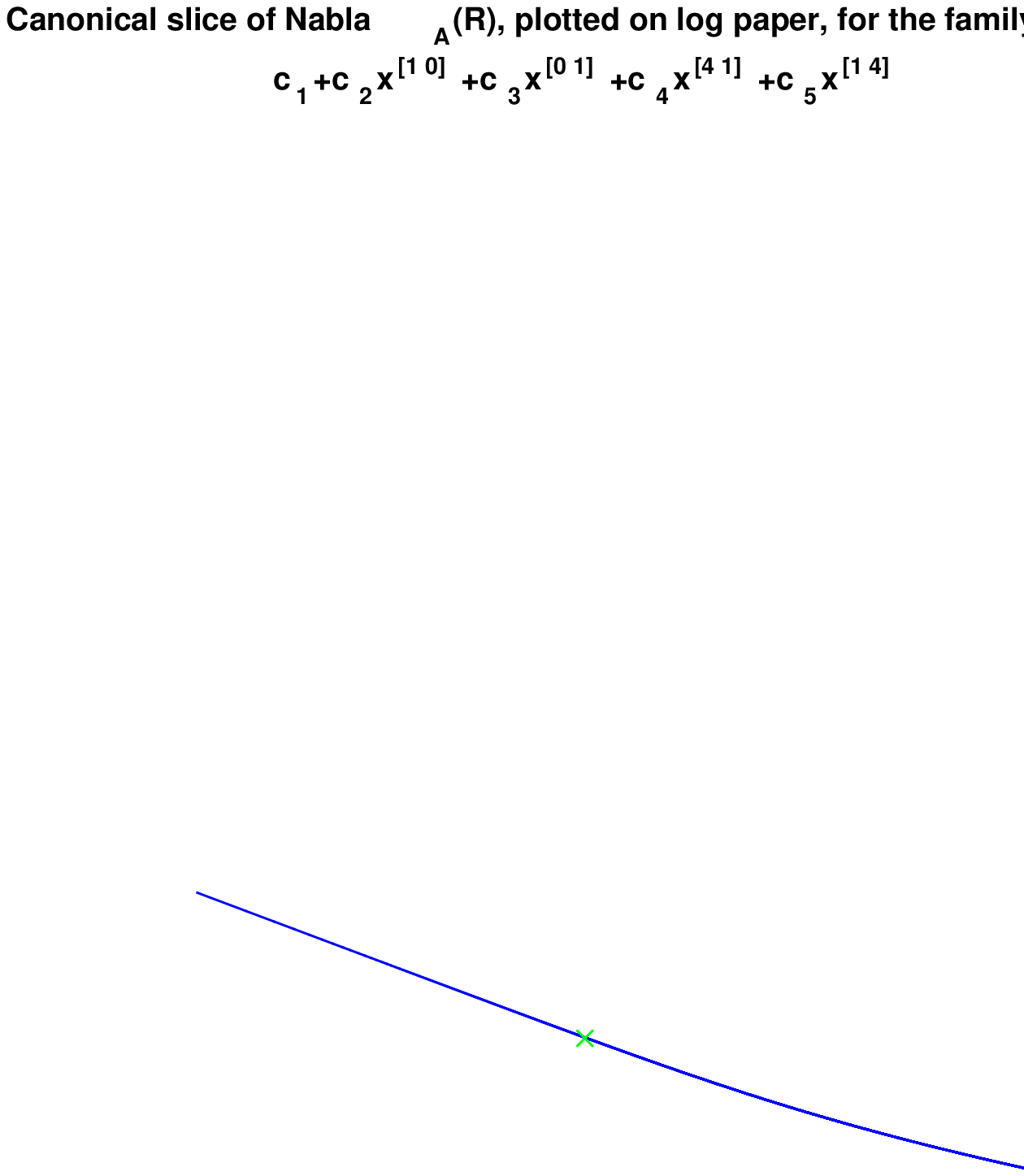,height=1.2in,clip=}}
\fbox{\epsfig{file=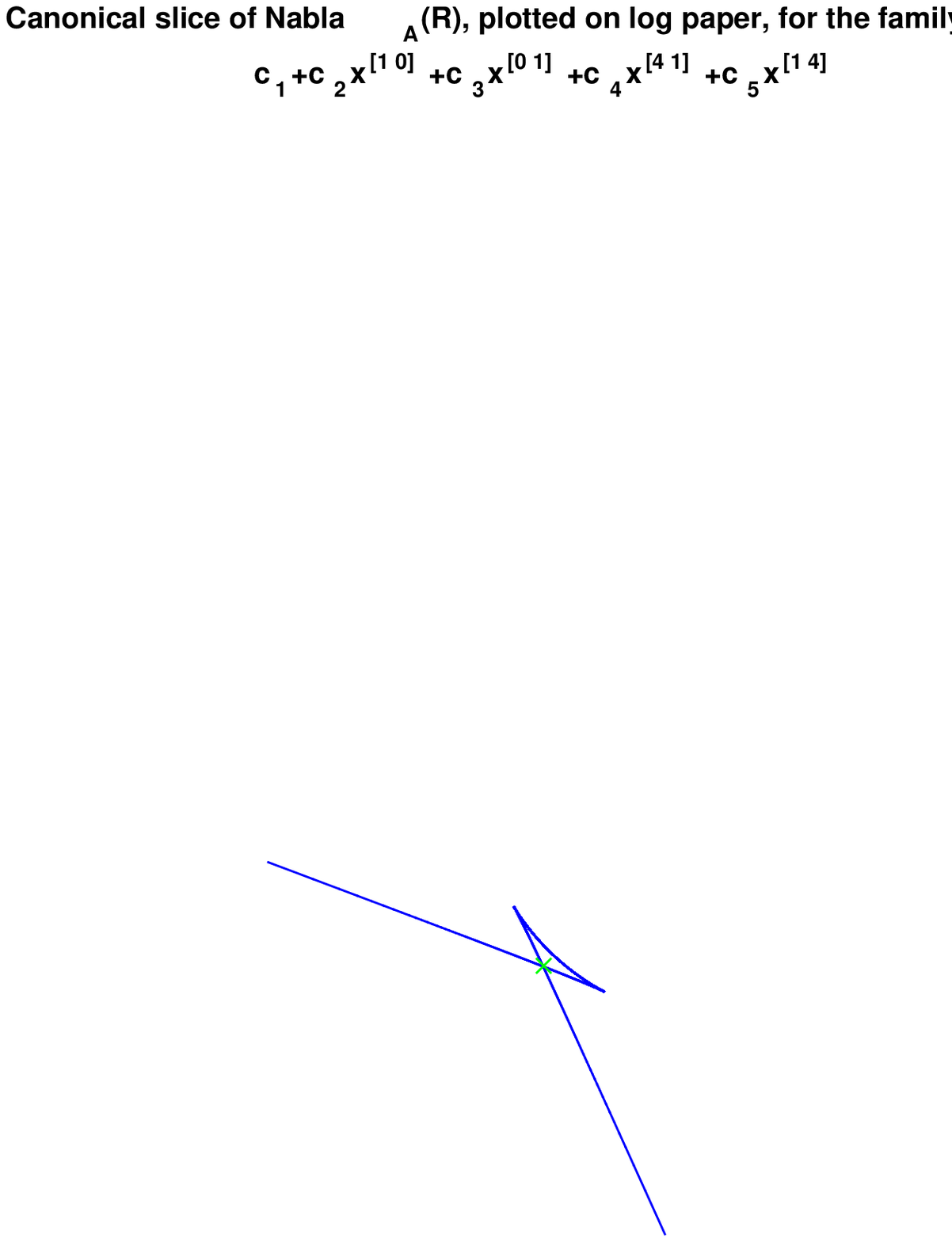,height=1.2in,clip=}}
\fbox{\epsfig{file=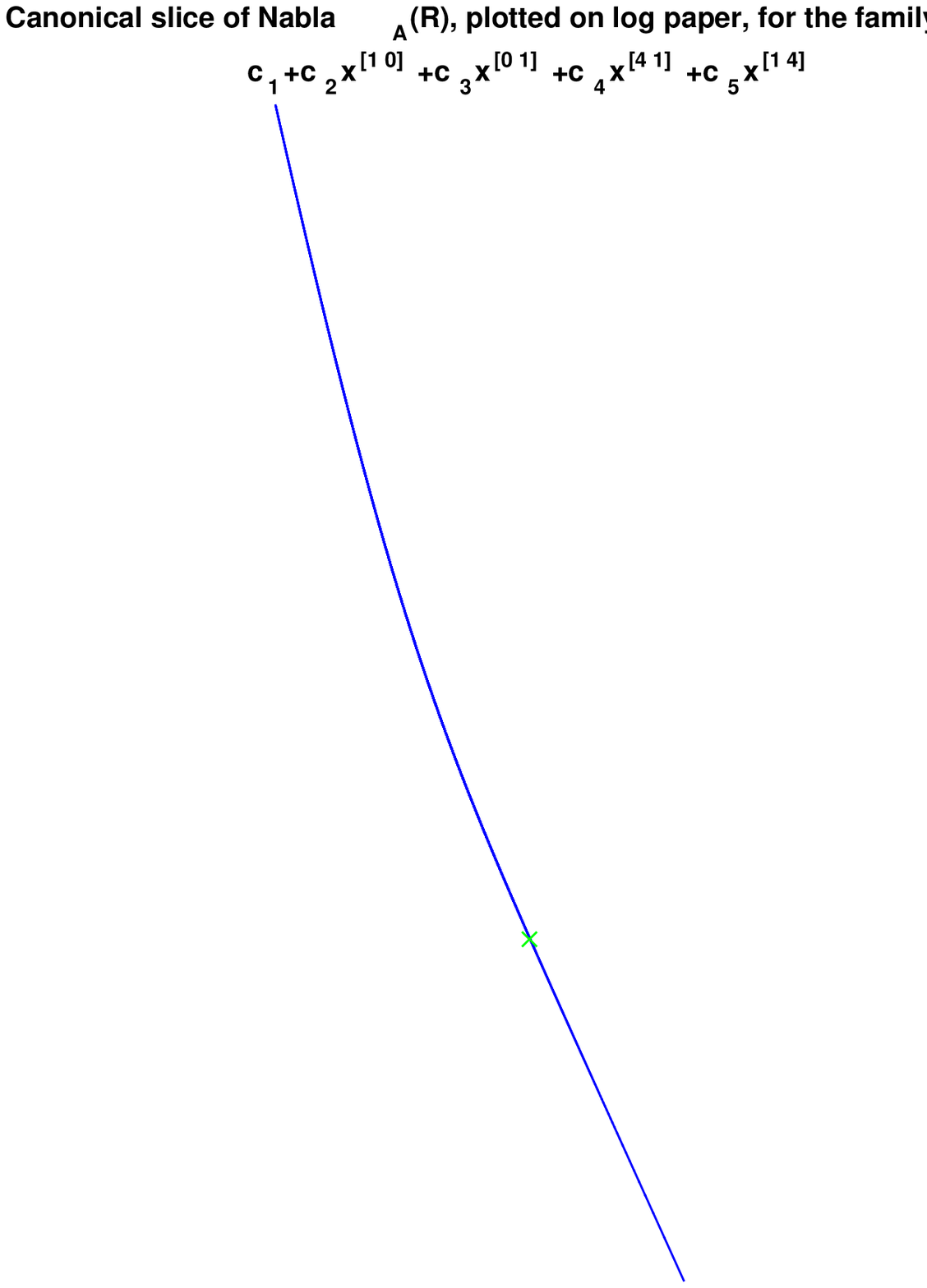,height=1.2in,clip=}}
\fbox{\epsfig{file=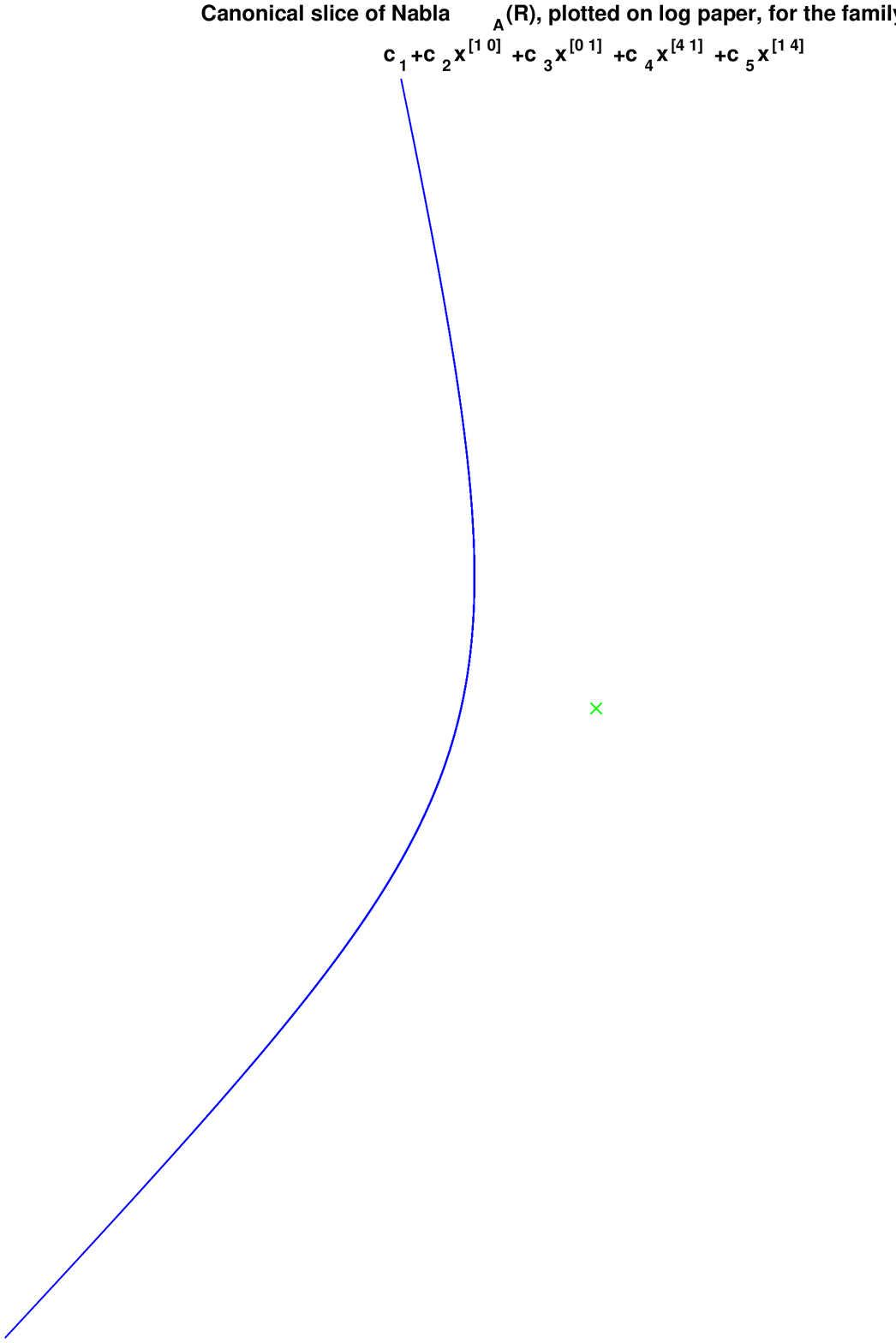,height=1.2in,clip=}}}
\put( 30,60){$++--+$}
\put(128,60){$-+++-$}
\put(230,-10){$+--++$}
\put(327,-10){$-++-+$}
\put(455,-10){$+-++-$}
\end{picture}}

\noindent 
Note that the curves drawn above are in fact unbounded, 
so the number of reduced signed chambers for the $\sigma$ above, from left to 
right, is respectively $2$, $2$, $3$, $2$, and $2$. (The tiny $\times$ in 
each illustration indicates the origin in $\R^2$.) In particular, 
only $\sigma\!=\!(1,-1,-1,1,1)$ yields an inner chamber. Note also that 
$\Gamma(\cA,B)$ is always the union of all the $\Gamma_\sigma(\cA,B)$.  
\begin{rem} While the shape of the reduced signed chambers 
certainly depends on the choice of $B$, the hyperplane arrangement 
$H_\cA$ and the number of signed chambers 
for any fixed $\sigma$ are independent of $B$. In particular, working 
with the $\Gamma_\sigma(\cA,B)$ in $\R^{n+k-d(\cA)-1}$ helps us visualize and 
work with $\Xi_\cA$, which lives in $\Pro^{n+k-1}_\R$. \dia 
\end{rem} 

\section{Morse Theory, Fewnomial Bounds, and the Proof of 
Theorem \ref{thm:big}}  
Let us call $\cA\!\in\!\R^{n\times (n+k)}$ {\em combinatorially simplicial}
if and only if $\cA\cap Q$ has cardinality $1+\dim Q$ for every face $Q$ of
$\conv \{a_1,\ldots,a_{n+k}\}$. 
(The books \cite{grunbaum,ziegler} are excellent standard references 
on polytopes, their faces, and their normal vectors.)  
Note that $\conv \{a_1,\ldots,a_{n+k}\}$ 
need {\em not} be a simplex for $\cA$ to be combinatorially simplicial 
(consider, e.g., Example \ref{ex:penta}). We now state the main reason we care 
about reduced signed chambers.  
\begin{thm}
\label{thm:morse} \cite{rojasrusek} 
Suppose $\cA\!\in\!\R^{n\times (n+k)}$ is combinatorially simplicial, 
non-defective, and $g_1$ and $g_2$ are each $n$-variate exponential 
$(n+k)$-sums with spectrum $\cA$. Suppose further that 
$\sign(c_{g_1})\!=\!\pm \sign(c_{g_2})$, and $(\Log|c_{g_1}|)B$ and 
$(\Log|c_{g_2}|)B$ lie in the same reduced discriminant chamber. Then 
$Z_\R(g_1)$ and $Z_\R(g_2)$ are ambiently isotopic in $\Rn$. \qed  
\end{thm} 

\noindent 
The special case $\cA\!\in\!\Z^{n\times (n+k)}$, without the use of 
$\Log$ or $B$, is alluded to near the beginning of 
\cite[Ch.\ 11, Sec.\ 5]{gkz94}. 
However, Theorem \ref{thm:morse} is really just an instance of 
{\em Morse Theory} \cite{milnor,smt},    
once one considers the manifolds defined by the fibers of the 
map\linebreak $Z_\R(g)\mapsto (\Log|c_g|)B$ along paths inside a fixed 
signed chamber. In particular, the assumption that $\cA$ be combinatorially 
simplicial forces any topological change in $Z_\R(g)$ to arise solely from 
singularities of $Z_\R(g)$ in $\Rn$. When $\cA$ is more general, topological 
changes in $Z_\R(g)$ can arise from pieces of $Z_\R(g)$ approaching infinity, 
with no singularity appearing in $\Rn$. So our chambers will need to be cut 
into smaller pieces. 

So we now address arbitrary $\cA$, but 
we'll first need a little more terminology. 
\begin{dfn} 
\label{dfn:chamber2}
Given any $\cA\!\in\!\R^{n\times (n+k)}$ with distinct 
columns, and any outer normal $w\!\in\!\Rn$ to a face of 
$\conv \cA$, we let $\cA^w\!:=\![a_{j_1},\ldots,a_{j_r}]$ denote the 
sub-matrix of $\cA$ corresponding to the set 
$\{a\!\in\!\cA \; | \; a\cdot w\!=\!\max_{a'\in\cA}\{a'\cdot w\}\}$.  
We call $\cA^w$ a {\em (proper) non-simplicial face} of $\cA$ when 
$d(\cA^w)\!\leq\!d(\cA)-1$ and $\cA^w$ has at least $d(\cA^w)+1$ columns. 
Also let $B^w$ be any matrix whose columns form a basis for the right 
nullspace of $\widehat{\left(A^w\right)}$, and let $\pi_w : \C^{n+k} 
\longrightarrow \C^r$ be the natural coordinate projection map defined by 
$\pi_w(c_1,\ldots,c_{n+k})\!:=\!(c_{j_1},\ldots,c_{j_r})$. When $\cA$ is 
non-defective and {\em not} combinatorially simplicial we then define the 
{\em completed} reduced signed contour,  
$\widetilde{\Gamma}_\sigma(\cA,B)\subset\R^{n+k-d(\cA)-1}$, 
to be the union of $\Gamma_\sigma(\cA,B)$ and\\
\mbox{}\hfill $\displaystyle{\bigcup\limits_{\substack{\cA^w \text{a non-}\\ 
\text{simplicial}\\ \text{face of } \cA}} 
\overline{\left\{\left.\pi^{-1}_w\!\left(\Log|\lambda(B^w)^\top|\right)B 
\; \right| \; \mathrm{sign}\left(\Log\left|\lambda (B^w)^\top\right|\right) 
\!=\!\pm \pi_w(\sigma)  \ , \ [\lambda]\!\in\!
\Pro^{n+k-d(\cA)-2}_\R\setminus H_\cA\right\}}}$.\hfill\mbox{}\\ 
We call any unbounded connected component of $\R^{n+k-d(\cA)-1}\setminus
\widetilde{\Gamma}_\sigma(\cA,B)$ an {\em outer chamber}. 
Finally, we define $\widetilde{\Gamma}(\cA,B)\!:=\!\bigcup\limits_{\sigma\in
\{\pm 1\}^{n+k}}\widetilde{\Gamma}_\sigma(\cA,B)$. \dia 
\end{dfn} 
\begin{ex}
\label{ex:inf}
When $\cA\!=$\scalebox{.6}[.6]{$\begin{bmatrix}0 & 1 & 0 & 2 & 0\\ 
0 & 0 & 1 & 0 & 2 \end{bmatrix}$} it is easy to find a $B$ yielding the 
following reduced contour $\Gamma(\cA,B)$ and completed reduced contour 
$\widetilde{\Gamma}(\cA,B)$:\\ 
\mbox{}\hfill\epsfig{file=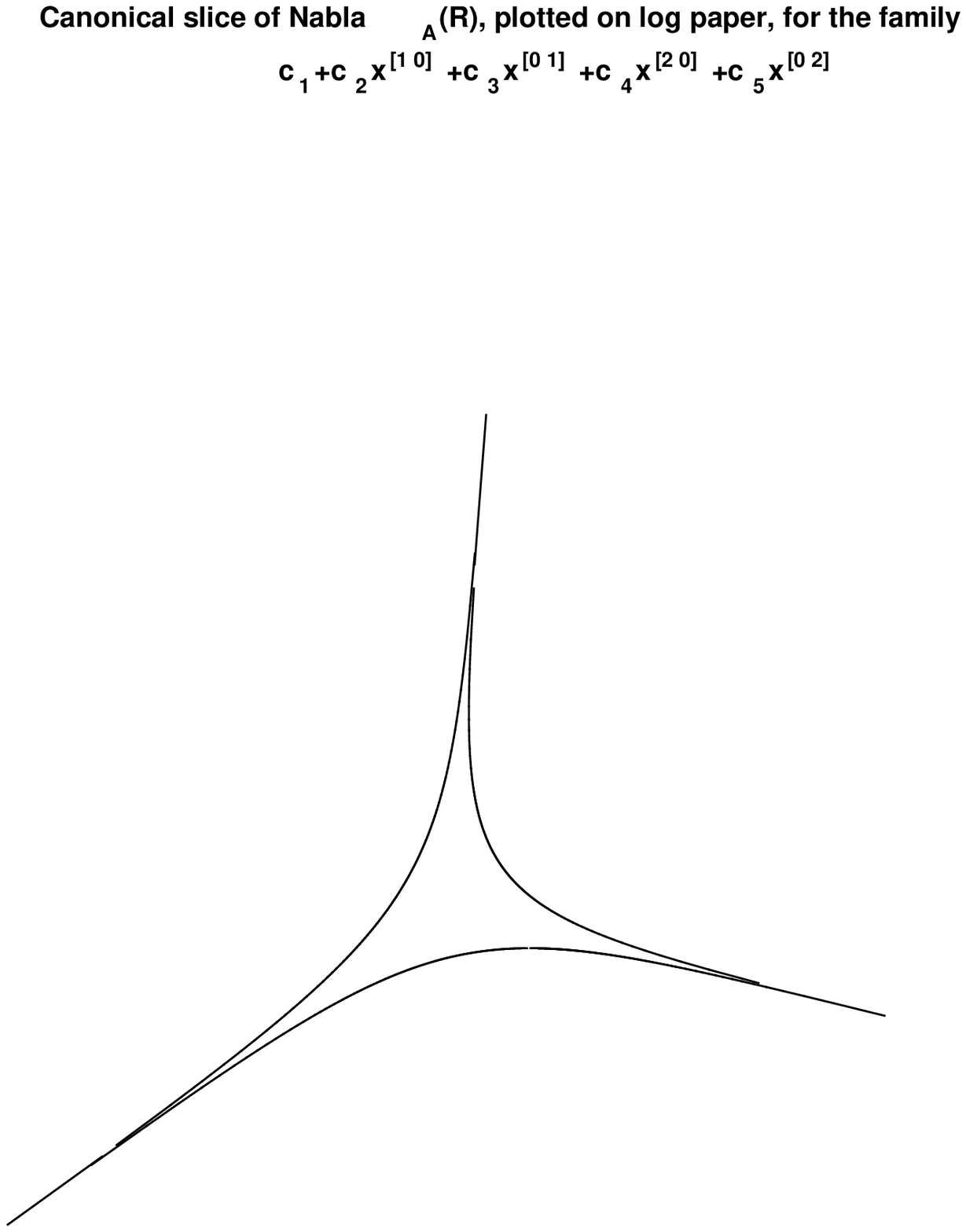,height=1.6in,clip=}
\hfill\epsfig{file=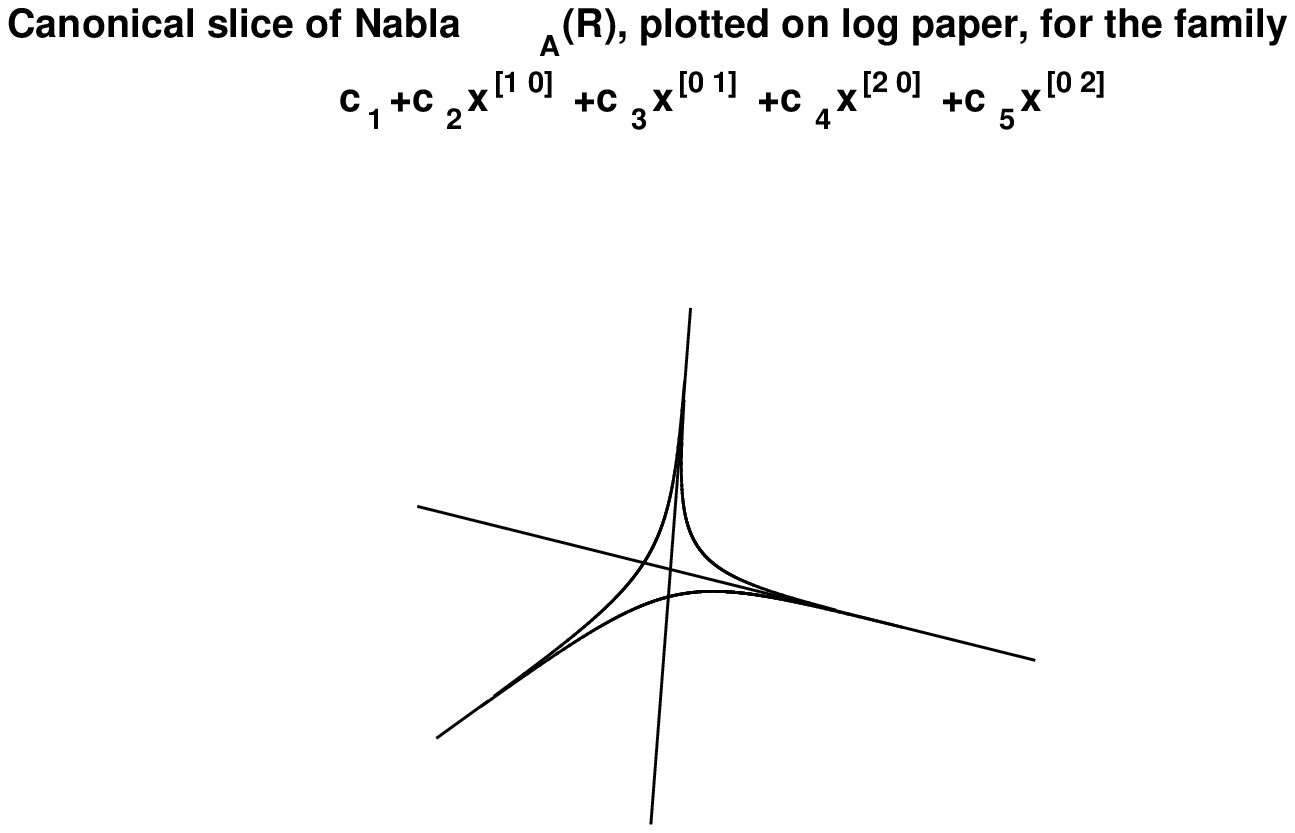,height=1.6in,clip=}
\hfill\mbox{}\\ 
Note in particular that $\widetilde{\Gamma}(\cA,B)\!=\!\Gamma(\cA,B)\cup 
S_1\cup S_2$ where $S_1$ and $S_2$ are lines that can be viewed as 
line bundles over points. These points are in fact $(\Log|\Xi_{\cA_1}|)B$ and 
$(\Log|\Xi_{\cA_2}|)B$ where $\cA_1$ and $\cA_2$ are the facets of $\cA$ with 
respective outer normals $(-1,0)$ and $(0,-1)$, and 
$B\!\approx$\scalebox{.6}[.6]
{$\begin{bmatrix}[r] 0.4335 & -0.8035 & -0.0635 & 0.4018 & 0.0317\\ 
0.3127 & 0.2002 & -0.8256 & -0.1001 & 0.4128\end{bmatrix}^\top$}. \dia
\end{ex}
\begin{prop}
\label{prop:non}
If $\cA$ is not combinatorially simplicial then
$\cA$ has at most $n+k-d(\cA)-1$ non-simplicial faces. \qed
\end{prop}
\begin{prop}
\label{prop:2face} 
\scalebox{.95}[1]{Suppose $k\!=\!3$, $\cA$ has exactly $2$ non-simplicial 
facets, $d(\cA)\!=\!n$, $B\!\in\!\R^{(n+3)\times 2}$}\linebreak 
is any matrix whose columns form a basis for the right 
nullspace of $\hA$, and $[\beta_{i,1},\beta_{i,2}]$ is\linebreak  
the $i\thth$ row of $B$. 
Then $\{[\beta_{i,1}:\beta_{i,2}]\}_{i\in\{1,\ldots,n+3\}}$ has cardinality 
$n+1$ as a subset of $\Pro^1_\R$, and $\widetilde{\Gamma}(\cA,B)
\setminus\Gamma(\cA,B)$ is a union of $2$ lines. \qed  
\end{prop} 
\begin{thm} \cite{rojasrusek} 
\label{thm:morse2} 
Suppose $\cA\!\in\!\R^{n\times (n+k)}$ is non-defective, {\em not} 
combinatorially simplicial, and $d(\cA)\!=\!n$. Suppose also that  
$g_1$ and $g_2$ are each $n$-variate exponential $(n+k)$-sums with 
spectrum $\cA$, $\sigma\!:=\!\sign(c_{g_1})\!=\!\pm \sign(c_{g_2})$,  
and $(\Log|c_{g_1}|)B$ and $(\Log|c_{g_2}|)B$  
lie in the same connected component of 
$\R^{n+k-d(\cA)-1}\setminus\widetilde{\Gamma}_\sigma(\cA,B)$. Then 
$Z_\R(g_1)$ and $Z_\R(g_2)$ are ambiently isotopic in $\Rn$. \qed   
\end{thm} 
\begin{ex} 
Observe that the circle defined by 
$\left(u+\frac{1}{2}\right)^2+\left(v-2\right)^2\!=\!1$ intersects the 
positive orthant, while the circle defined by 
$\left(u+\frac{3}{2}\right)^2+\left(v-\frac{3}{2}\right)^2=1$ does not. 
Consider then $\cA\!=$\scalebox{.6}[.6]{$\begin{bmatrix}0 & 1 & 0 & 2 & 0\\ 
0 & 0 & 1 & 0 & 2 \end{bmatrix}$} as in our last example, and let
\textcolor{blue}{$g_1\!=\!\left(e^{y_1}+\frac{1}{2}\right)^2
+\left(e^{y_2}-2\right)^2-1$} and\linebreak  
\textcolor{red}{$g_2\!=\!\left(e^{y_1}+\frac{3}{2}\right)^2+\left(e^{y_2}
-\frac{3}{2}\right)^2-1$}. Then \textcolor{blue}{$g_1$} and \textcolor{red}
{$g_2$} have spectrum $\cA$, $\sign(g_1)\!=\!\sign(g_2)\!=\!\sigma$ with 
$\sigma\!=\!(1,1,-1,1,1)$, and \textcolor{blue}{$(\Log|c_{g_1}|)B$} and 
\textcolor{red}{$(\Log|c_{g_2}|)B$} lie in the same reduced signed 
$\cA$-discriminant chamber (since 
$\Gamma_\sigma(\cA,B)\!=\!\emptyset$ here). However, \textcolor{blue}
{$Z_\R(g_1)$} consists of a single smooth arc, while 
\textcolor{red}{$Z_\R(g_2)$} is empty. 
This is easily explained by the {\em completed} contour 
$\widetilde{\Gamma}_\sigma(\cA,B)$ consisting of two lines, and 
\textcolor{blue}{$(\Log|c_{g_1}|)B$} and \textcolor{red}{$(\Log|c_{g_2}|)B$} 
lying in distinct connected components of 
$\R^2\setminus\widetilde{\Gamma}_\sigma(\cA,B)$ 
as shown, respectively via the symbols $\circ$ and $*$, below to 
the right. \dia 
\end{ex} 

\vspace{-.2cm} 
\noindent
\begin{minipage}[t]{.7 \textwidth}
\vspace{0pt}
Although we defined signed contours via a transcendental 
parametrization, they obey certain tameness properties akin to algebraic sets. 
One fundamental result implying this tameness is the following refined 
fewnomial bound. 
\end{minipage}\hspace{1cm} 
\begin{minipage}[t]{.2 \textwidth}
\vspace{0pt}
\begin{picture}(0,0)(0,0)
\put(-25,-60){\epsfig{file=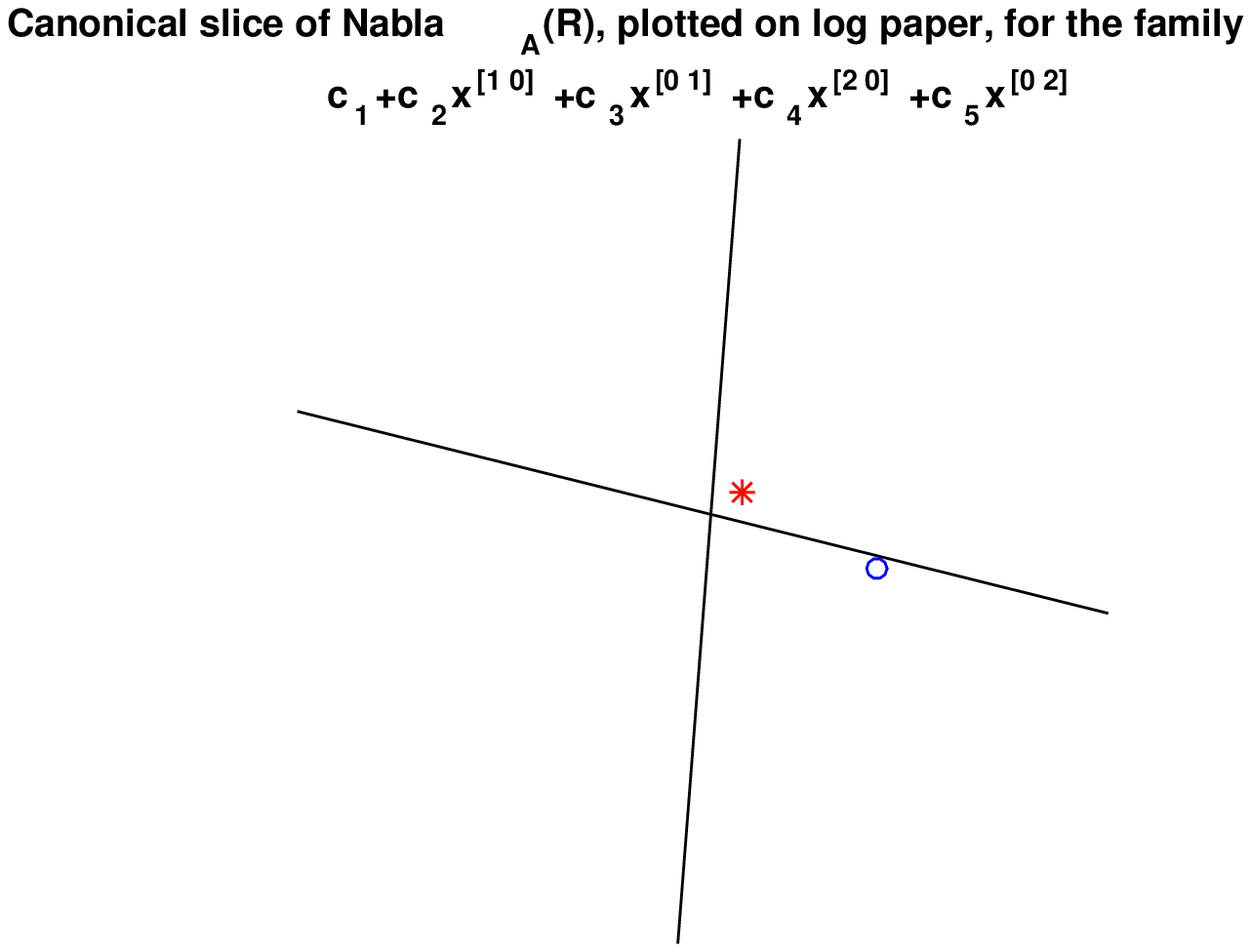,height=.9in,clip=}} 
\end{picture} 
\end{minipage} 

\begin{thm} \label{thm:gale} (See \cite[Thm.\ 3.1]{bihansottilegale},  
\cite{bbs}, \& \cite[Lem.\ 1.8]{barbados}) 
Suppose $m\!\geq\!1$,\linebreak 
$j\!\geq\!2$, $E\!=\![e_{i,\ell}]\!\in\!\R^{(m+j)\times j}$, 
$U:=[u_{i,\ell}]\!\in\!\R^{(m+j)\times (j+1)}$ has 
i$\thth$ row $(u_{i,0},u_{i,1},\ldots,u_{i,j})$, $u_i\!:=\!(u_{i,1},
\ldots,u_{i,j})$, $\Delta\!:=\!\left\{y\!\in\!\R^j\; | \; 
u_{i,0}+u_i\cdot y\!>\!0 \text{ for all } i\!\in\!\{1,\ldots,j\}\right\}$,
and \\ 
\mbox{}\hfill 
$\displaystyle{H\!:=\!\left(\prod^{m+j}_{\ell=1} (u_{\ell,0}+u_\ell\cdot y 
)^{e_{\ell,1}}, \ldots,\prod^{m+j}_{\ell=1} (u_{\ell,0}+u_\ell\cdot y 
)^{e_{\ell,j}}\right)-(1,\ldots,1)}$.\hfill\mbox{}\\  
Then $H$ has fewer than 
$S(m,j)\!:=\!\frac{e^2+3}{4}\left(\sqrt{2^{j-1}}m\right)^j$  
non-degenerate roots in $\Delta$. Furthermore, for $j\!=\!1$, 
$H$ has at most $S(m,1)\!:=\!m+1$ non-degenerate roots in $\Delta$, and there 
exist $H$ attaining $m+1$ distinct roots in $\Delta$. 
\qed 
\end{thm}  

\noindent 
We call systems of the above form {\em $j$-variate Gale Dual systems with 
$m+j$ factors}. 

\begin{cor} 
\label{cor:line} Suppose $\cA\!\in\!\R^{n\times (n+k)}$ is combinatorially 
simplicial, non-defective, $d(\cA)\!=\!n$, and $\sigma\!\in\!\{\pm 1\}^{n+k}$. 
Then, following the notation of Theorem \ref{thm:gale}, a generic affine line 
$L\subset\R^{k-1}$ intersects  
$\Gamma_\sigma(\cA,B)$ in no more than $\left\lfloor S(n+2,k-2)
\right\rfloor$ points when $k\!\geq\!3$. Also, for $k\!=\!2$ there 
is at most $S(n+2,0)\!:=\!1$ intersection. 
\end{cor} 

\noindent 
{\bf Proof of Corollary \ref{cor:line}:} When $k\!=\!2$ 
we have that $\Gamma(\cA,B)$ is merely a point, so this case 
follows easily. So let us assume $k\!\geq\!3$ and let 
$[\cL_{i,j}]_{(i,j)\in\{1,\ldots,k-2\}\times \{0,\ldots,k-1\}}
\!\in\!\R^{(k-2)\times k}$ be any matrix defining the affine line $L$ as 
follows:\\  
\mbox{}\hfill $L=\left\{x\!\in\!\R^{k-1}\; | \;
\cL_{i,1}x_1+\cdots+\cL_{i,k-1}x_{k-1}\!=\!\cL_{i,0} \text{ for all }
i\!\in\!\{1,\ldots,k-2\}\right\}$.\hfill\mbox{}\\ 
Also let $(\xi_1,\ldots,\xi_{k-1})\!:=\!\xi_{\cA,B}$. 
(So each $\xi_i$ is a logarithm of the absolute value of a 
linear form in $\lambda_1,\ldots,\lambda_{k-1}$.) 
Note then that $L$ meets $\Gamma(\cA,B)$ at the point 
$\xi_{\cA,B}([\lambda])$ only if 
\begin{eqnarray*} 
\sum^{k-1}_{\ell=1} \cL_{1,\ell}\xi_\ell([\lambda]) & = & \cL_{1,0}\\ 
 & \vdots & \\ 
\sum^{k-1}_{\ell=1} \cL_{k-2,\ell}\xi_\ell([\lambda]) & = & \cL_{k-2,0}
\end{eqnarray*} 

\noindent 
Exponentiating both sides of the preceding system, and collecting
factors, we obtain that\linebreak 
\scalebox{.97}[1]{there is a matrix
$E\!=\![E_{i,j}]\!\in\!\R^{(k-2)\times (n+k)}$ such that $L$ meets
$\Gamma(\cA,B)$ at the point $\xi_{\cA,B}([\lambda])$ only if} 
\begin{eqnarray*} 
\prod^{n+k}_{\ell=1} 
  (\beta_\ell\cdot \lambda)^{E_{1,\ell}}&=& e^{\cL_{1,0}}\\ 
 & \vdots & \\ 
\prod^{n+k}_{\ell=1} 
   (\beta_\ell\cdot \lambda)^{E_{k-2,\ell}}&=& e^{\cL_{k-2,0}}
\end{eqnarray*} 

\noindent 
Setting $\lambda_{k-1}\!=\!1$ to dehomogenize the linear forms $\beta_i
\cdot \lambda$, Theorem \ref{thm:gale} then tells us that $L$ meets 
$\Gamma(\cA,B)$ at no more than $S(n+2,k-2)$ points. Since the number of 
intersections is an integer, we can take floor and conclude. \qed 

\begin{lemma} 
\label{lemma:sums} 
If $n,k',k''\!\geq\!2$ then $S(n+1,k'+k''-2)\!\geq\!S(n+1,k'-2)+S(n+1,k''-2)$. 
More generally, if $k_1+\cdots+k_r\!=\!k-1$ with 
$k_i\!\geq\!2$ for all $i$ and $r\!\geq\!2$, then 
$S(n+1,k-5)+1\!\geq\!\sum^r_{i=1}S(n+1,k_i-2)$. 
\end{lemma} 

\noindent 
\scalebox{.93}[1]{{\bf Proof of Lemma \ref{lemma:sums}:} The first 
assertion is immediate since $S(n+1,k'-2)+S(n+1,k''-2)$}\linebreak 
$\leq\!2S(n+1,k''-2)$ (assuming $k''\!\geq\!k'$) and 
$2^{1+(k''-2)(k''-3)/2}\!\leq\!2^{(k'+k''-2)(k'+k''-3)/2}$. 
The second assertion follows easily by induction: Writing 
$k\!=\!(\cdots ((k_1+k_2)+k_3)+ \cdots+k_{r-1})+k_r$, the first 
assertion of our lemma implies that $\sum^r_{i=1}S(n+1,k_i-2)\!\leq
\!S(n+1,k'-2)+S(n+1,k''-2)$ for some $k',k''\!\geq\!2$ 
with $k-1\!=\!k'+k''$. It is then easy to see (from the 
power of $2$ factor of $S(m,j)$ again) that 
$S(n+1,k'-2)+S(n+1,k''-2)\!\leq\!S(n+1,k-3-2)+S(n+1,2-2)$, 
i.e., the left-hand side of the inequality is maximized 
when $\{k',k''\}\!=\!\{2,k-3\}$. \qed 

\begin{cor}
\label{cor:line2} Suppose $\cA\!\in\!\R^{n\times (n+k)}$ is 
non-defective, $\cA$ is {\em not} combinatorially simplicial, 
$d(\cA)\!=\!n$, and $\sigma\!\in\!\{\pm 1\}^{n+k}$. Then a generic affine 
line $L\subset\R^{k-1}$ intersects $\widetilde{\Gamma}_\sigma(\cA,B)$ in no 
more than\\  
\scalebox{.89}[1]{$S(n+2,k-2)+S(n+1,k-5)+\cdots+S\!\left(n+2-\min\left\{n+1,
\left\lfloor\frac{k-2}{3}\right\rfloor\right\},k-2-3
\min\left\{n+1,\left\lfloor\frac{k-2}{3}\right\rfloor\right\}\right)$}  
\linebreak 
\mbox{}\hspace{1cm}$+\min\left\{n+1,\left\lfloor\frac{k-2}{3}\right\rfloor
\right\}$\\  
points when $k\!\geq\!4$. Also, for $k\!\in\!\{2,3\}$ we have respective upper 
bounds of $1$ and $n+5$. 
\end{cor}

\noindent
{\bf Proof of Corollary \ref{cor:line2}:} We simply follow 
essentially the same argument as the proof of Corollary \ref{cor:line}, 
save that we work with $\widetilde{\Gamma}_\sigma(\cA,B)$ instead of 
$\Gamma_\sigma(\cA,B)$. In particular, the case $k\!=\!2$ presents no new 
difficulties since $\widetilde{\Gamma}_\sigma(\cA,B)$ is always a point. The 
case $k\!=\!3$ follows easily upon observing, thanks to Proposition 
\ref{prop:2face}, that $\widetilde{\Gamma}_\sigma(\cA,B)\setminus\Gamma_\sigma
(\cA,B)$ is either empty, a line, or two lines. 

For $k\!\geq\!4$ we simply observe that $L$ will either 
intersect $\Gamma_\sigma(\cA,B)$ or some fiber closure of the form 
$\overline{\left\{\left.\pi^{-1}_w\!\left(\Log|\lambda(B^w)^\top|\right)B 
\; \right| \; \mathrm{sign}\left(\Log\left|\lambda (B^w)^\top\right|\right) 
\!=\!\pm \pi_w(\sigma)  \ , \ [\lambda]\!\in\!
\Pro^{n+k-d(\cA)-2}_\R\setminus H_\cA\right\}}$. There are no more than 
$S(n+2,k-2)$ of the first kind of intersection, thanks to Corollary 
\ref{cor:line}. After applying the map $\pi_w$, we see that counting the 
second kind of intersections reduces to a lower-dimensional instance of 
Corollary \ref{cor:line}. In particular, the second kind of intersections, 
for fixed $w$, contribute no more than $S(\dim(\cA^w)+2,k_w-2)$ to our total, 
where $k_w$ is the number of columns of $\cA^w$ minus $d(\cA^w)$. 
Note that the sum of all the $k_w$ is no more than $k-1$ since 
$d(\cA)\!=\!n$. Note also that when $\cA$ has just 
two non-simplicial facets, with one having exactly $n+1$ columns, 
the other has at most $n+k-4$ colums. In which case, 
these facets would contribute $S(n+1,0)+S(n+1,k-5)$ to our sum. 
In particular, this is the maximal possible contribution, over all 
distributions of points to the non-simplicial facets, thanks to Lemma 
\ref{lemma:sums}.  

More generally, the non-simplicial faces of $\cA$ naturally form a poset under 
containment which, along with the distribution of the columns of $\cA$ 
as points in the relative interior of the faces of $\conv\{a_1,
\ldots,a_{n+k}\}$, determines the sum of $S(m,j)$ giving an 
upper bound for the intersection count we seek. 
Lemma \ref{lemma:sums} then tells us that our sum is maximized when 
it is of the form\\ 
\mbox{}$S(n+2,k-2)+(S(n+1,0)+S(n+1,k-5))+\cdots$\\ 
\scalebox{.87}[1]{$\cdots+\left(S\!\left(n+2
 -\min\left\{n+1,\left\lfloor\frac{k-2}{3}\right\rfloor\right\},
 0\right)+S\!\left(n+2-\min\left\{n+1,\left\lfloor\frac{k-2}{3}\right\rfloor
 \right\},k-2-3\min\left\{n+1,\left\lfloor\frac{k-2}{3}\right\rfloor\right\}
 \right)\right)$.}\\ 
Since $S(m,0)\!=\!1$ for all $m$ we are done. \qed 

\begin{thm}
\label{thm:jens}
\cite{jens} 
Let $[c_g]$ be any smooth point of $\Xi_\cA$. Then 
$Z_\R(g)$ has a unique singular point $\zeta$,  
and the Hessian of $g$ at $\zeta$ has full rank. \qed 
\end{thm} 

In what follows, let $N(g)$ denote the number of connected 
components of $Z_\R(g)$. 
\begin{thm}
\label{thm:exptrop} \cite{rojasrusek} 
If $g$ is an $n$-variate exponential $(n+k)$-sum with 
spectrum $\cA\!\in\!\R^{n\times (n+k)}$, and $(\Log|c_g|)B$ lies in an outer 
chamber, then $N(g)\!\leq\!(n+k)(n+k-1)/2$. \qed 
\end{thm} 
\begin{thm} 
\label{thm:side} 
Suppose $n\!\geq\!2$ and $g_-,g_*,g_+$ are $n$-variate exponential $(n+k)$-sums 
with non-defective spectrum $\cA$, 
$\sign(c_{g_-})\!=\!\sign(c_{g_*})\!=\!\sign(c_{g_+})\!=\!\sigma$,  
and $L'\subset\R^{k-1}$ is the unique line segment connecting 
$(\Log|c_{g_-}|)B$ and $(\Log|c_{g_+}|)B$. Suppose further that 
$L'\cap\widetilde{\Gamma}_\sigma(\cA,B)\!=\!\{(\Log|c_{g_*}|)B\}$, 
and $(\Log|c_{g_*}|)B$ is a smooth point of 
$\widetilde{\Gamma}_\sigma(\cA,B)$. Then $|N(g_+)-N(g_-)|\!\leq\!1$ and 
$|N(g_\pm)-N(g_*)|\!\leq\!1$.  
\end{thm} 

\noindent 
{\bf Proof:}  
$X\!:=\!\{(c_g,y)\!\in\!\R^{n+k}\times \Rn\; 
| \; g(y)\!=\!0 \ , \ \sign(c_g)\!=\!\sigma \ , \ (\Log|c_g|)B\!\in\!L'\}$ 
forms a singular real manifold but, thanks to Theorem \ref{thm:jens}, 
$X$ has a unique singularity at $(c_{g_*},\zeta)$ where $\zeta\!\in\!\Rn$ 
is the unique singular point of $Z_\R(g_*)$. Let $\phi : [-1,1]\longrightarrow 
\R^{n+k}$ be any smooth function with $\sign(\phi(t))\!=\!\sigma$ 
for all $t\!\in\![-1,1]$ and $(\Log|\phi([-1,1])|)B\!=\!L'$. 
Let $\pi : \R^{n+k}\times \Rn \longrightarrow \R^{n+k}$ denote the natural 
orthogonal projection forgetting the second factor.  
We then see that $\phi^{-1}\circ \pi$ is a Morse function on $X$. 
By Stratified Morse Theory \cite{milnor,smt}, there is a closed ball 
$U\!\subset\!\R^{n+k}\times \Rn$ 
containing $(c_{g_*},\zeta)$ such that $U\cap X$ is homeomorphic to a real 
hypersurface of the form $Y\!=\!\{(x,t)\!\in\!\Rn\times [-1,1]\; | 
\; Q(x)\!=\!t, |x|\!\leq\!1\}$, where $Q$ is a homogeneous quadratic form 
with signature identical to the Hessian of $g_*$ at $\zeta$, 
$Y\cap\{t\!=\!\pm 1\}$ is isotopic to $U\cap Z_\R(g_\pm)$, and 
$Y\cap\{t\!=\!0\}$ is isotopic to $U\cap Z_\R(g_*)$. 

To conclude, observe that $Y\cap \{t\!=\!\pm 1\}$ empty implies 
that the signature of $Q$ is $\pm (1,\ldots,1)$, and thus 
$Y\cap\{t\!=\!0\}$ is a point and $Y\cap \{t\!=\!\mp 1\}$ is a 
sphere. So then $U\cap Z_\R(g_\pm)$ empty implies that 
$U\cap Z_\R(g_\mp)$ has a unique isolated connected component. 
In other words, the conclusion of our theorem is true. 

If $Y\cap \{t\!=\!\pm 1\}$ are both non-empty, then the signature 
of $Q$ can {\em not} be $\pm (1,\ldots,1)$. So then 
$Y\cap \{t\!=\!-1\}$, 
$Y\cap \{t\!=\!0\}$, and 
$Y\cap \{t\!=\!1\}$, each have at least one connected component, 
and none has more than $2$ connected components. This in turn 
implies that $U\cap Z_\R(g_-)$, $U\cap Z_\R(g_*)$, and  
$U\cap Z_\R(g_+)$ each have at least one connected component, 
and none has more than $2$ connected components. Note also that 
any connected component of $U\cap Z_\R(g_\pm)$ (resp.\linebreak 
\scalebox{.96}[1]{$U\cap Z_\R(g_*)$) 
lies in a unique connected component of $Z_\R(g_\pm)$ (resp.\ $Z_\R(g_*)$). 
So we are done. \qed}   

\subsection{The Proof of Theorem \ref{thm:big}:}  
\label{sub:proof} 
If $n\!=\!1$ then the theorem follows easily from the well-known generalization 
of Descartes' Rule of Signs to real exponents (see, e.g., \cite{wang}), 
and with an improved (tight) upper bound of $k$. So let us assume henceforth 
that $n\!\geq\!2$. 

\medskip 
\noindent 
{\bf Combinatorially Simplicial Case:} 
If $\cA$ is defective then $\Xi_\cA\cap \Pro^{n+k-1}_\R$ has real 
codimension $2$ in $\Pro^{n+k-1}_\R$ and thus $\Pro^{n+k-1}_\R\setminus 
\Xi_\cA$ is path-connected. So then, by the framework of our proof of 
Theorem \ref{thm:side}, the number of connected components of $g$ is constant 
for any fixed choice of sign vector. So it suffices to count connected 
components in outer chambers and, by Theorem \ref{thm:exptrop}, we are done. 
So let us now assume $\cA$ is non-defective.  

Consider a line segment $L_{gh}$, connecting $(\Log|c_g|)B$ to $(\Log|c_h|)B$, 
where $h$ has the same spectrum as $g$ and $\sign(c_h)\!=\!\sign(c_g)\!=:\!
\sigma$, but known cardinality for $Z_\R(h)$. The key trick will then be that 
$L_{gh}$ intersects $\Gamma_\sigma(\cA,B)$ in few places, and the number of 
connected components of an $f$ with $\Log|c_f|\!\in\!L$ changes only slightly 
as $f$ moves from $h$ to $g$. 

In particular, we may assume in addition that $h$ lies in an 
outer chamber $\cC_\sigma$ (since outer chambers are open and unbounded). 
By Theorem \ref{thm:getchambers} we may then assume that $L_{gh}$ lies in an 
affine line $L$ sufficiently generic for Corollary \ref{cor:line} to hold,  
{\em and} that $L_{gh}$ intersects $\Gamma_\sigma(\cA,B)$ only at smooth points 
of $\Gamma_\sigma(\cA,B)$. Furthermore, since the points of 
$L_{gh}\cap \Gamma_\sigma(\cA,B)$ can be linearly 
ordered, we may also assume that $(h,\cC_\sigma)$ has been 
chosen so that $L_{gh}\cap \Gamma_\sigma(\cA,B)$ consists of 
no more than half of $L\cap \Gamma_\sigma(\cA,B)$. 

If we can show that $Z_\R(h)$ has few connected components, and $Z_\R(f)$ 
gains few connected components as $f$ moves from $h$ to $g$ (with 
$\Log|c_f|$ restricted to $L$), then we'll be done. 

Toward this end, observe that $Z_\R(h)$ has at most 
$(n+k)(n+k-1)/2$ connected components, thanks to Theorem \ref{thm:exptrop}. 
Since we have chosen $L_{gh}$ so that it intersects $\Gamma_\sigma(\cA,B)$ 
only at smooth points, Theorem \ref{thm:side} tells us that as $f$ moves from 
$h$ to $g$ (with $(\Log|c_f|)B$ restricted to $L$), each such intersection 
introduces at most $1$ new connected component. (Theorems \ref{thm:morse}  
also tell us that $N(f)$ is constant when $(\Log|c_f|)B$ lies between adjacent 
intersections in $L\cap \Gamma_\sigma(\cA,B)$.) So by Corollary \ref{cor:line}, 
we are done with the case where $\cA$ is combinatorially\linebreak   
\scalebox{.96}[1]{simplicial, with a slightly smaller upper bound of 
$\displaystyle{\frac{(n+k)(n+k-1)}{2}+\left\lfloor S(n+2,k-2)/2\right\rfloor}$. 
\qed}  

\medskip 
\noindent 
{\bf The Case Where $\cA$ is not Combinatorially Simplicial:} 
Here we just slightly modify the argument we used when $\cA$ was 
combinatorially simplicial: The key difference is that we 
work with $\widetilde{\Gamma}_\sigma(\cA,B)$ instead of $\Gamma_\sigma(\cA,B)$, 
and apply Corollary \ref{cor:line2} instead of Corollary \ref{cor:line}.  
The number of intersections $L$ with 
$\widetilde{\Gamma}_\sigma(\cA,B)$ between $(\Log|c_g|)B$ and 
$(\Log|c_h|)B$ then clearly admits an upper bound of\\  
$T(n,k)\!:=\!(S(n+2,k-2)+S(n+1,k-5)+$\\ 
\mbox{}\hspace{1.8cm}\scalebox{.88}[1]{$\cdots+S\!\left(n+2-\min\left\{n+1,
\left\lfloor\frac{k-2}{3}\right\rfloor\right\},k-2-3
\min\left\{n+1,\left\lfloor\frac{k-2}{3}\right\rfloor\right\}\right)
+\min\left\{n+1,\left\lfloor\frac{k-2}{3}\right\rfloor\right\})/2$.}  

At this point, we are nearly done, but for some elementary observations on  
sums of powers of $2$ and the size of $S(n+1,1)$. First, observe that the 
powers of $2$ in the summands making up $T(n,k)$ are:  
\[ 2^{(k-2)(k-3)/2},2^{(k-5)(k-6)/2},\ldots,2^{
\left.\left(k-2-3\min\left\{n+1,\left\lfloor\frac{k-2}{3}\right\rfloor\right\}
\right)
\left(k-3-3\min\left\{n+1,\left\lfloor\frac{k-2}{3}\right\rfloor\right\}
\right)\right/2
}.\]  
So, in particular, the sum of all but the first power of $2$ is 
strictly less than
\[ 2^{k^2-11k+30}+2^{k^2-11k+29}+\cdots+2^4+2^3
=2^{k^2-11k+31}-8\!<\!2^{(k-4)(k-5)/2}-8.\]  
Next, we observe that $\min\left\{n+1,\left\lfloor\frac{k-2}{3}\right\rfloor
\right\}\!\leq\!n+1\!<\!S(n+1,1)$. So then we easily obtain that 
$T(n,k)\!\leq\!(S(n+2,k-2)+S(n+1,k-4))/2$. So the final upper 
bound we obtain is $N(g)\!\leq\!\frac{(n+k)(n+k-1)}{2}+\left\lfloor 
(S(n+2,k-2)+S(n+1,k-4))/2\right\rfloor$, which is slightly better than 
our stated bound. \qed 

\section{Acknowledgements} 
We humbly thank Erin Lipman for asking an insightful question that 
was the genesis for this paper: ``What bounds are known for chamber-depth, as 
opposed to the number of chambers?'' We also thank Saugata Basu and Frank 
Sottile for valuable discussions on genericity.

\end{document}